\documentclass[12pt]{article}
\usepackage{latexsym}
\usepackage{theorem}
\usepackage{graphicx}
\usepackage{amsmath,color}
\usepackage{amsfonts}
\usepackage{natbib}
\usepackage{soul}

\theorembodyfont{\rmfamily}
\newtheorem{theorem}{Theorem}

\newtheorem{lemma}[theorem]{Lemma}

\newtheorem{definition}[theorem]{Definition}

\theoremstyle{break}

\newtheorem{example}[theorem]{Example}

\title{Searching for Multiple Objects in Multiple Locations}

\author{Thomas Lidbetter\thanks{Department of Management Science and Information Systems, Rutgers Business School, Newark, NJ 07102, tlidbetter@business.rutgers.edu} \and Kyle Y. Lin\thanks{Operations Research Department, Naval Postgraduate School, Monterey, CA 93943, kylin@nps.edu}}

\begin{document}

\maketitle

\begin{abstract}
\noindent
Many practical search problems concern the search for multiple hidden objects or agents, such as earthquake survivors.
In such problems, knowing only the list of possible locations, the Searcher needs to find all the hidden objects by visiting these locations one by one.
To study this problem, we formulate new game-theoretic models of discrete search between a Hider and a Searcher.
The Hider hides $k$ balls in $n$ boxes, and the Searcher opens the boxes one by one with the aim of finding all the balls.
Every time the Searcher opens a box she must pay its {\em search cost}, and she either finds one of the balls it contains or learns that it is empty.
If the Hider is an adversary, an appropriate payoff function may be the expected total search cost paid to find all the balls, while if the Hider is Nature, a more appropriate payoff function may be the difference between the total amount paid and the amount the Searcher would have to pay if she knew the locations of the balls {\em a priori} (the {\em regret}).
We give a full solution to the regret version of this game, and a partial solution to the search cost version. We also consider variations on these games for which the Hider can hide at most one ball in each box. The search cost version of this game has already been solved in previous work, and we give a partial solution in the regret version.
\end{abstract}



\newpage

\section{Introduction}
Consider search problems in which a {\em Searcher} must find $k \ge 1$ balls distributed among a set of $n > k$ boxes. Each box has a {\em search cost}, which is the cost the Searcher must pay to open the box.
If the Searcher opens an empty box, she finds no ball and learns the box is empty.
If the Searcher opens a nonempty box, then she finds just one ball without learning whether the box contains any more.
The objective is to find a randomized search to minimize the worst-case expected cost the Searcher has to pay to find all $k$ balls. Equivalently, this problem can be formulated as a zero-sum game between the Searcher and a malevolent cost-maximizing {\em Hider}.

It is natural in many situations that the number of objects found depends on the amount of effort invested in searching. For example, consider the search for a row of landmines in a particular site.
Since minefields often exhibit regular linear patterns with equal spacing \citep{Lake-Keenan}, the landmines are found one by one until reaching the end of the row, whereupon the Searcher learns that there are no more hidden. 
Another example is after an earthquake, when the survivors may be hidden in several collapsed buildings, and the search and rescue team needs to dig deeper in the rubble to save them one by one.
A further example, which comes from the natural world, is that of animals or birds that hide, or {\em cache}, food for later rediscovery \citep{VanderWall}. Squirrels that bury nuts in the ground (the Hider) risk those nuts later being pilfered by a competitor (the Searcher).

When the Hider is an adversary---such as an enemy force who plants landmines---it is reasonable to adopt a payoff function equal to the expected total cost to find all the balls.
The expected total cost, however, need not be the only interesting payoff function to study, especially when the Hider is Nature.
Consider an example with $n=2$ boxes and $k=5$ balls, and the search cost is \$10 for box 1 and \$1 for box 2.
The Hider can easily guarantee a total search cost at least \$50 simply by hiding all the 5 balls in box 1.
But if the Hider is Nature, then the Searcher need not feel particularly unhappy if she finds all 5 balls by opening box 1 for 5 straight times without opening box 2 at all---because the Searcher did not \textit{make a mistake}.
In these situations, a more appropriate objective function is to minimize the \textit{regret}---the difference between the search cost and the cost she would paid to find the balls if she already knew their locations.
In economics, there is extensive study on how people make choices based on regret~\citep{kahneman1979prospect, loomes1982regret, bell1982regret}.

Since the Searcher needs to open the same box several times to find all the balls hidden there, we call our model a \textit{multi-look} search game.
In contrast, a variation to our model is one in which the Hider can put at most one ball in each box, so that the Searcher does not need to look in any box more than once. We call this a \textit{single-look} search game.
By considering the two different rules (\textit{multi-look} or \textit{single-look}) and the two different payoff functions (\textit{search cost} or \textit{regret}), we have a total of four search games.
In this paper, we study three of these four search games, while the single-look search game with search cost has been previously solved in ~\cite{lidbetter}.
Table~\ref{tab:results} summarizes the known results for these four versions of the game.

\begin{table}[htb]
	\caption{Summary of results.}
	\begin{center}
		\begin{tabular}{|l|c|c|}
			\hline
			& {\bf Search Cost} & {\bf Regret} \\ \hline
			{\bf Multi-look}  & Solved for $n=2$ and for equal costs & Full solution \\
		 & (Section~\ref{sec:cost-multi}) 	&  (Section~\ref{sec:regret-multi})  \\  
			\hline
			{\bf Single-look} & Full solution & Solved for $k=n-1$    \\
	& (\cite{lidbetter}) 	&  (Section~\ref{sec:regret-single})  \\ 
			\hline
		\end{tabular}
		\label{tab:results}
	\end{center}
\end{table}


Our work falls into the category of {\em search games}, on which there is a rich literature; see, for example, \cite{AlpernGal}, \cite{Gal2011} 
and \cite{Hohzaki} for surveys. Although many search games are played in continuous time and space, discrete search games are also well studied. In particular, there has been a lot of work on searching for a {\em single} object in a finite set of hiding locations. In a classic problem solved by Blackwell (reported in \cite{Matula}), the objective is to minimize the expected cost of finding an object hidden in boxes according to a known probability distribution, where each box has a search cost and an overlook probability.
A game-theoretic version of this game was considered in \cite{Bram}, \cite{Roberts:1978wl}, \cite{Gittins:1979vj}, \cite{Ruckle}, and more recently in \cite{Lin-Singham}.

Search problems in which there are several hidden objects are less well studied. \cite{Assaf} and \cite{Sharlin} study a variation of Blackwell's problem (described above) in which the Searcher has to find {\em one} of several hidden objects in least expected cost. \cite{lidbetter} solves a game in which all the objects hidden in a set of boxes must be found in minimum expected cost. {\em Caching games} were introduced in \cite{caching}, in which a Searcher with limited resources aims to maximize the probability of locating a certain number of hidden objects. Caching games share a feature of our multi-look games, in that the number of objects found in a location can depend on how much searching effort is invested.

In some of the games studied in this paper, the payoff function we use is the \textit{additive regret} incurred by the Searcher. There is a lot of work in search theory that concentrates on the {\em multiplicative regret}---the \textit{ratio} of the cost paid by the Searcher in finding a hidden item to the cost she would have paid if she already knew its location. This type of payoff function was introduced in search theory independently by \cite{beck} and \cite{bellman} in the context of the {\em linear search problem}, later extended by \cite{gal74}. The linear search problem is framed in a continuous setting, but multiplicative regret for search theory has also been studied in a discrete setting in \cite{fixed} and \cite{searchratio}.

The rest of this paper proceeds as follows.
Sections~\ref{sec:cost-multi}, \ref{sec:regret-multi}, and \ref{sec:regret-single} is each concerned with a version of the search game as summarized in Table~\ref{tab:results}.
Whereas we solve the multi-search game with regret in Section~\ref{sec:regret-multi}, for the other two versions we can only derive the optimal solution in a few cases, and show by examples the complexity of the optimal solution in general.
Section~\ref{sec:conclusion} concludes and points out a few future research directions.

\section{Multi-look Search Game with Search Cost}
\label{sec:cost-multi}
Consider the multi-look search game with search cost, as summarized in Table~\ref{tab:results}. The Hider distributes $k$ balls among $n$ boxes in any way he wants, including putting up to $k$ balls in the same box.
The Searcher opens the boxes one by one, each time paying a search cost.
For conciseness, we write ${[n]:=\{1,\ldots,n\}}$ for the set of boxes, and let $c_i$ denote the search cost of box $i \in [n]$.
If the box contains at least one ball, the Searcher finds one of them but does not learn whether there are any others in the box; otherwise, she learns the box is empty. The payoff is the sum of the search cost paid after the Searcher finds all $k$ balls.
The Searcher is the minimizer and the Hider is the maximizer.

Since multiple balls can be hidden in a box, a pure strategy for the Hider is a nonnegative integer solution to the equation $\sum_{i=1}^n x_i = k$, where $x_i$ corresponds to the number of balls hidden in box $i \in [n]$.
Hence, the Hider has ${n+k -1 \choose k}$ pure strategies.
The Searcher's strategy set is harder to describe. 
Since the Searcher opens the boxes one by one, she can use the search history to decide which box to open next.
The search history at each stage can be described by the sequence of boxes that have been opened up to that point and those search outcomes.
As an example, suppose $n=k=2$.
A possible Searcher strategy could be described as follows.
Open box 1 first; if it contains a ball then open it again, and then open box 2 if the balls have not both already been found; if box 1 does not contain a ball then open box 2 twice.
Since the strategy set of each player is finite, the game has a value and optimal strategies.
Since the Searcher's strategy set is much larger than the Hider's, the Searcher may have many strategies that achieve optimality.
As will be seen in this paper, the Searcher often has an optimal strategy that does not depend on the full search history, and will be simpler to describe.





Let $T_k = T_k([n])$ denote the sum over all products $\Pi_{i=1}^n c_i^{x_i}$, for all nonnegative integer vectors $(x_1, \ldots, x_n)$ summing to $k$.
That is,
\[
T_k([n]) = \sum_{x_1+\cdots+x_n =k} \left( \prod_{i=1}^n c_i^{x_i} \right).
\]
For instance,
\[
T_2([3]) = c_1^2 + c_2^2 + c_3^2 + c_1 c_2 + c_2 c_3 + c_1 c_3.
\]

Consider the Hider strategy, which chooses a pure strategy $\mathbf{x} = (x_1, \ldots, x_n)$ with probability
\begin{equation}
p_{n, k}(\mathbf{x}) =  \frac{\prod_{i=1}^n c_i^{x_i}}{T_k([n])}.
\label{eq:HiderEqual}
\end{equation}
It is quick to check that this mixed strategy has an important property: whether or not the Searcher finds a ball after opening a box, the other balls are still hidden according to (\ref{eq:HiderEqual}).
It turns out that the mixed strategy in (\ref{eq:HiderEqual}) is an \textit{equalizing strategy}---a strategy that results in the same payoff against any Searcher strategy.
Note that if $c_i=c$ for all $i \in [n]$, this strategy simply chooses uniformly between all the Hider's pure strategies.

\begin{lemma}
	\label{le:Unk}
	For the muli-look search game with expected cost, against the Hider strategy $p_{n,k}$ in (\ref{eq:HiderEqual}), any Searcher strategy has expected cost
	\begin{equation}
	U([n],k) \equiv \frac{k T_{k+1}([n])}{T_k([n])}.
	\label{eq:Unk}
	\end{equation}
\end{lemma}
\textit{Proof.}
The proof is by induction on $n+k$. For $n=k=1$, we have 
\[
\frac{1 \cdot T_{2}([1])}{T_1([1])} =  \frac{1 \cdot c_1^2}{c_1}= c_1,
\]
which is clearly the expected cost, since the Searcher will open box 1 exactly once to find the only ball, and be done with it.

Suppose $n+k  \geq 3$ and suppose without loss of generality that the Searcher starts by opening box $n$, which costs $c_n$.
As noted already, whether or not a ball is found, the remaining balls will be hidden according to $\mathbf{p}$ in (\ref{eq:HiderEqual}).
More precisely, if a ball is found (which happens with probability $c_n T_{k-1}([n])/T_k([n])$, the remaining balls are hidden according to $p_{n,k-1}$; if a ball is not found (which happens with probability $T_k([n-1])/T_k([n])$, the remaining balls are hidden according to $p_{n-1,k}$. 
Hence,
\begin{align*}
U([n],k) & = c_n + \frac{c_n T_{k-1}([n])}{T_k([n])} U([n],k-1) + \frac{T_k([n-1])} {T_k([n])}  U([n-1],k)  \\
& = c_n + \frac{c_n T_{k-1}([n])} {T_k([n])} \left( \frac{(k-1) T_{k}([n])}{T_{k-1}([n])}  \right) + \frac{T_k([n-1])}{T_k([n])} \left( \frac{k T_{k+1}([n-1])}{T_k([n-1])} \right) \\
& \qquad (\mbox{by induction})&\\
& = \frac{ c_n T_{k}([n]) + c_n (k-1) T_k([n]) + k T_{k+1}([n-1])}{T_k([n])} \\
& = \frac{ k (c_n T_k([n]) + T_{k+1}([n-1]))}{T_k([n])} \\
& =  \frac{k T_{k+1}([n])}{T_k([n])},
\end{align*}
which completes the proof.
\hfill $\Box$

\smallskip

A natural conjecture is that the equalizing strategy in (\ref{eq:HiderEqual}) is optimal for the Hider. Indeed, this conjecture is true when the search costs are equal for all boxes, as we will show in Subsection~\ref{sec:equal}. However, the equalizing strategy is not optimal in general, as we show in Subsection~\ref{sec:n=2}.

\subsection{Special Case: Equal Search Costs} \label{sec:equal}
This section presents the solution to the special case of the game where all the search costs are equal.
Note that the single-look version of the game with equal search costs is trivial, as it has a simple solution where both players simply randomize uniformly between all their strategies.

To simplify notation, we normalize $c_i=1$ for $i \in [n]$.
In this game, $T_k([n]) = {n+k-1 \choose n-1}$, so by (\ref{eq:Unk}), the Hider can ensure that the value $V(n,k)$ of the game satisfies
\begin{align}
V(n,k) &\ge \frac{k{n+k \choose n-1}}{{n+k-1 \choose n-1}} = (n+k) \left( 1 - \frac{1}{k+1} \right) . \label{eq:equal-costs}
\end{align}
We will show that the preceding is indeed the value of the game, so that the Hider's equalizing strategy in~(\ref{eq:HiderEqual}) of choosing uniformly between all his pure strategies is optimal. The Searcher's optimal strategy is less obvious. A first guess at an optimal Searcher strategy would be to choose at each stage uniformly at random between all remaining boxes that may contain a ball. However, this strategy is not optimal, as we illustrate in the following example.

Suppose $n=k=2$. The Hider ensures a payoff of $(n+k)(1-1/(k+1)) = 8/3$ by using his equalizing strategy.
Consider the Searcher strategy that always chooses a random box that may contain a ball.
If the 2 balls are in the first box, then the Searcher has to pay a cost of 2 for opening the first box twice.
In addition, the probability she also opens the second box at some stage is $1-(1/2)(1/2) = 3/4$. Hence, the expected payoff is $2+3/4 > 8/3$.

In this example, an optimal strategy for the Searcher is as follows.
Open a box at random. If it does not contain a ball then open the other box twice.
If the box does contain a ball, then with probability $2/3$ open it again then open the other box if both balls have not yet been found; with probability $1/3$ open the other box then open the first box again if both balls have not yet been found. If the 2 balls are in the same box, this strategy has expected payoff $(1/2)(3)+(1/2)(1+(2/3)\cdot 1 + (1/3) \cdot 2) = 8/3$. If the balls are in different boxes, the expected payoff is $1+(2/3)\cdot 2 + (1/3) \cdot 1 = 8/3$.
Theorem~\ref{thm:eqcost} describes the optimal Searcher strategy for arbitrary $n$ and $k$.

\begin{theorem} \label{thm:eqcost}
	For the multi-look search game with equal search costs $c_i = 1$, the optimal strategy for the Hider is given by (\ref{eq:HiderEqual}). The optimal Searcher strategy is recursive. If there are $k$ balls still to be found and $n$ boxes left that may contain balls, the Searcher chooses a box uniformly at random, and with probability $p_j$ she opens that box until either finding that it is empty or finding $j$ balls (whichever happens sooner), where $p_j = \lambda(k+1-j)$, $j=1,\ldots,n$ and $\lambda = 2/(k(k+1))$ is a normalizing factor. The value of the game is 
	\[
	V(n,k) = (n+k) \left( 1 - \frac{1}{k+1} \right).
	\]
\end{theorem}
\noindent
\textit{Proof.}
Write $\Gamma(n,k)$ for this multi-look search game with $n$ boxes and $k$ balls.
By (\ref{eq:equal-costs}), it is sufficient to show that the Searcher strategy described in the statement of the theorem ensures a payoff of at most $(n+k)/(1 - (k+1))$. We prove the result by induction on $n+k$, noting that the case with $n=k=1$ is trivial. 

Consider an arbitrary pure strategy for the Hider, for which there are $n_i$ boxes containing $i$ balls, $i=0,\ldots,k$. Note that we must have
\begin{align}
\sum_{i=0}^{k}n_i &= n, \label{eq:sum1}\\
\sum_{i=0}^{k}i n_i & = k. \label{eq:sum2}
\end{align}
With probability $n_i/n$, the Searcher will pick a box containing $n_i$ balls. In this case, for $j \le i$, with probability $p_j$, the Searcher will pay $j$, find $j$ balls, and will then be faced with the game $\Gamma(n,k-j)$, which, by induction, has value $V(n,k-j) = (n+k-j)/(1 - 1/(k-j+1))$. For $j > i$, with probability $p_j$, the Searcher will pay $i+1$, find $k$ balls, and will then by faced with the game $\Gamma(n-1,k-i)$, which, by induction, has value $V(n-1,k-i) = (n+k-i-1)/(1 - 1/(k-i+1))$. 

Putting this together, we obtain
\begin{align*}
V(n,k) & \le \sum_{i=0}^k \frac{n_i}{n} \left( \sum_{j=1}^i p_j (j + V(n,k-j)) + \sum_{j=i+1}^k p_j (i+1 + V(n-1,k-i)) \right) \\
& = \frac{\lambda}{n} \sum_{i=0}^k n_i \left( \sum_{j=1}^i (k+1-j)\left(n+k - \frac{n+k-j}{k-j+1}\right)  + \sum_{j=i+1}^k (k+1-j) \left( n+k - \frac{n+k-i-1}{k-i+1} \right) \right) \\
& =\frac{\lambda(n+k)}{n} \sum_{i=0}^k n_i \sum_{j=1}^k (k+1-j) - \frac{\lambda}{n} \sum_{i=0}^k n_i \left( \sum_{j=1}^i (n+k-j) + \left( \frac{n+k-i-1}{k-i+1} \right) \sum_{j=i+1}^k k+1 - j \right) \\
& = n+k - \frac{\lambda}{2n} \sum_{i=0}^k n_i(k^2 + kn -k + in) \qquad \text{ (by (\ref{eq:sum1}))}\\
& = n+k - \frac{\lambda}{2n} ((k^2 + kn - k)n + nk) \qquad \text{ (by (\ref{eq:sum1}) and (\ref{eq:sum2}))}\\ 
& = (n+k) \left( 1 - \frac{1}{k+1} \right),
\end{align*}
which completes the proof.
\hfill $\Box$



\subsection{Special Case: $n=2$} \label{sec:n=2}
We now relax the assumption that all the search costs are the same, and we will see that the Hider's equalizing strategy is not necessarily optimal, even for $n=2$. Note that in this case, when a search is unsuccessful, the Searcher learns all the remaining balls are hidden in the other box, so without writing it explicitly, we can assume that after the Searcher finds a box is empty, she will search in the other box until finding all the remaining balls.

Suppose the Searcher chooses between $k+1$ pure strategies, where the pure strategy indexed by some $j \in \{0 ,1,\ldots, k\}$ starts by opening box 1 at most $j$ times.
If the Searcher finds $j$ balls in box 1, then she opens box 2 until finding all balls or finding box 2 empty, in which case she returns to box 1 to find all balls.
If the Searcher finds box 1 empty before finding $j$ balls, then she opens box 2 to find all the other balls.
If the Hider hides $i$ balls in box 1, and the Searcher opens box 1 at most $j$ times in her first $k$ searches, then the total cost to find all balls is
\[
C(i,j) = \left\{
\begin{array}{ll}
i  c_1 + (k-i) c_2 + c_1, & \mbox{ if } i < j, \\
i  c_1 + (k-i) c_2 , & \mbox{ if } i = j, \\
i  c_1 + (k-i) c_2 + c_2, & \mbox{ if } i > j. \\
\end{array}
\right.
\]

%

\begin{lemma}
	\label{le:CostEqualizing}
	Consider the muli-look search game with search cost with $n=2$ boxes. If
	\begin{equation}
	q_{k}(j) = 1 - \frac{k c_1^{k-j} c_2^{j}}{T_k([2])} \geq 0,
	\label{eq:qkj}
	\end{equation}
	for $j=0, \ldots, k$, then the Searcher's optimal mixed strategy is to open box 1 at most $j$ times with probability $q_k(j)$.
	The Hider's optimal mixed strategy is given in (\ref{eq:HiderEqual}), and the value of the game is
	\begin{equation}
	\frac{k T_{k+1}([2])}{T_k([2])}.
	\label{eq:n=2,value}
	\end{equation}
\end{lemma}
\textit{Proof.}
First, it is straightforward to check $\sum_{j=0}^k q_k(j) = 1$.
Next, if the Hider hides $i$ balls in box 1, then the Searcher's mixed strategy in (\ref{eq:qkj}) yields expected cost
\begin{align*}
\sum_{j=0}^k q_k(j) C(i,j) &= i c_1 + (k-i) c_2 + c_1 \sum_{j=i+1}^k q_k(j) + c_2 \sum_{j=0}^{i-1} q_k(j) \\
&= i c_1 + (k-i) c_2 + (k-i)c_1 - c_1 \sum_{j=i+1}^k \frac{k c_1^{k-j} c_2^{j}}{T_k([2])} + i c_2 - \sum_{j=0}^i \frac{k c_1^{k-j} c_2^{j}}{T_k([2])} \\
&= k(c_1 + c_2) - k \left( \frac{c_1 c_2 T_{k-1}([2])}{T_k{[2]}} \right)\\
&= \frac{k T_{k+1}([2])}{T_k([2])}.
\end{align*}
Since the Hider and the Searcher each has a mixed strategy that guarantees the same payoff in (\ref{eq:n=2,value}) regardless of what the other does, the proof is completed.
\hfill $\Box$

\smallskip

From now on we assume $c_1 \ge c_2$, without loss of generality.
Since $q_k(j)$ increases in $j$, it is a legitimate probability distribution if $q_k(0) \geq 0$, or equivalently,
\[
-(k-1) + r + r^2 + \cdots + r^k \geq 0.
\]
where $r = c_2/c_1$.

\begin{lemma}
	For $m \geq 2$, define
	\[
	f_m(r) \equiv -(m-1) + r + r^2 + \cdots + r^m.
	\]
	The function $f_m(r)$ has a unique root in $(0,1)$.
	Writing $r_m$ as the unique root in $(0,1)$ for $f_m(r)$, then the sequence $r_2, r_3, \ldots$ is increasing.
\end{lemma}
\textit{Proof.}
Since $f_m(r)$ is continuous and increases in $r$, with $f_m(0) = -(m-1) < 0$ and $f_m(1)=1> 0$, it follows that it has a unique root in $(0,1)$.
In addition, since
\[
f_{m+1}(r_m) = f_m(r_m) - 1 + r_m^{m+1} < 0,
\]
it follows that $r_{m+1} > r_m$, which completes the proof.
\hfill $\Box$

\smallskip


In the multi-look search game with expected cost and $n=2$, if $c_2/c_1 \geq r_k$, then $q_k(0) \geq 0$, so Lemma~\ref{le:CostEqualizing} applies and we have the optimal policy.
If $c_2/c_1 < r_k$, then the Hider can guarantee an expected payoff of at least $U([2],k)$ but the Searcher cannot guarantee the payoff does not exceed that using the strategy of Lemma~\ref{le:CostEqualizing}, which suggests that the value of the game is higher than that.

Indeed, consider the example $c_1=10$, $c_2=1$, and $k=2$. In this case, $U([2],2) = 2(1000+100+10+1)/(100+10+1)\approx 20.0$. Consider an alternative strategy for the Hider of placing one of the balls in box 1 with probability 1, and then choosing where to put the second ball by using his equalizing strategy~(\ref{eq:HiderEqual}) for $n=2,k=1$. Since a best response to this strategy would be to first open box 1 then open the two boxes in any order, this would guarantee an expected payoff of at least $10+U([2],1) = 10+(100+10+1)/(10+1) \approx 20.1$.
The fact that $10+U([2],1) > U([2], 2)$ suggests that in general the Searcher should hide a subset of the balls in box 1 with probability 1, if $c_1$ is sufficiently larger than $c_2$.



How many balls should the Hider set aside and put in box 1?
If the Hider sets aside $m$ balls in box 1, and then distributes the other $k-m$ balls between the 2 boxes according to his equalizing strategy~(\ref{eq:HiderEqual}), then the Hider guarantees expected cost $m c_1 + U([2], k-m)$.
It turns out that the best strategy of this type is the Hider's optimal strategy, and the value of the game is
\[
\max_{m=0, \ldots, k} m c_1 + U([2], k-m).
\]
The optimal number of balls to set aside in box 1 depends on $c_2/c_1$ through the function $f_m(c_2/c_1)$, as explained in the next theorem.


\begin{theorem}
	Consider the muli-look search game with expected cost where there are $n=2$ boxes with costs $c_1,c_2, c_1 \ge c_2$. Let $b =  \max \{m: r_m \leq c_2 /c_1 \}$.
	If $k \leq b$, then an optimal strategy for the Hider is given by~(\ref{eq:HiderEqual}), and an optimal strategy for the Searcher is given by~(\ref{eq:qkj}); the value of the game is $U([2], k)$.
	If $k > b$, then an optimal strategy for the Hider is to first set aside $k-b$ balls and put them in box 1, and then hide the other $b$ balls according to (\ref{eq:HiderEqual}).
	An optimal strategy for the Searcher is to first open box 1 for up to $k-b$ times.
	If any of these searches finds box 1 to be empty, then search in box 2 until finding the remaining balls.
	If the Searcher finds $k-b$ balls in box 1, then she uses the mixed strategy according to $q_b(j)$ in (\ref{eq:qkj}) for the rest of the search.
	The value of the game is $(k-b) c_1 + U([2], b)$.

\end{theorem}
\textit{Proof.}
If $k \leq b$, then $c_2/c_1 \geq r_k$ and $q_k(0) \geq 0$, so the result follows from Lemma~\ref{le:CostEqualizing}.
Now suppose $k > b$.
Since the stated Hider's strategy guarantees expected cost at least $(k-b) c_1 + U([2], b)$ (Lemma~\ref{le:Unk}), it remains to show that the stated Searcher's strategy guarantees an expected cost at most $(k-b) c_1 + U([2], b)$ for any Hider's strategy.

Suppose the Hider hides $i$ balls in box 1.
Consider two cases:
\begin{enumerate}
	\item
	$i \geq k-b$: The Searcher will first find $k-b$ balls in box 1, incurring $(k-b)c_1$ in the process.
	Thereafter, the Searcher uses the mixed strategy according to $q_b(j)$ in (\ref{eq:qkj}), which is a legitimate probability distribution, so the result follows from Lemma~\ref{le:CostEqualizing}.
	\item
	$i < k-b$:
	The Searcher will first find $i$ balls in box 1 in her first $i+1$ searches, incurring $(i+1) c_1$ in the process, and learns that all the other $k-i$ balls are in box 2.
	Therefore, the total cost is $(i+1)c_1 + (k-i)c_2$.
	To show that this cost is less than $(k-b) c_1 + U([2],b)$, we take the difference to get
	\begin{align*}
	(i+1)c_1 + (k-i)c_2 - (k-b) c_1 - U([2],b) &= (i+1-k+b)(c_1-c_2) + (b+1)c_2 - U([2],b) \\
	&\leq (b+1)c_2 - U([2],b),
	\end{align*}
	since $i+1 \leq k-b$ and $c_1 \ge c_2$.
	Writing $r=c_2/c_1$, the right-hand side of the preceding can be written as
	\begin{align*}
	c_1 \left( (b+1) r - \frac{b (1+r+\cdots+r^{b+1})}{1+r+\cdots+r^b} \right) &= \frac{c_1 f_{b+1}(r)}{1+r+\cdots+r^b}\\
	& < \frac{c_1 f_{b+1}(r_{b+1})}{1+r+\cdots+r^b}\\
	&=0,
	\end{align*}
	where the inequality follows from the definition of $r$ and $b$.
\end{enumerate}
Consequently, the stated Searcher strategy guarantees an expected cost at most $(k-b) c_1 + U([2], b)$ for any Hider strategy, which completes the proof.
\hfill $\Box$

\subsection{The General Case}

We might conjecture that in the general case, an optimal strategy for the Hider would involve setting aside some number of balls to be placed in the highest-cost boxes with probability 1, and then distributing the remaining balls among the boxes according to the equalizing strategy in~(\ref{eq:HiderEqual}).
This conjecture, however, is not true, as seen in the following example.

\begin{example}
\label{ex:CostEx1}
Consider the multi-look search game with search cost, where $n=4$, $k=2$, and $(c_1, c_2, c_3, c_4) = (10, 9, 1, 1)$.
Using linear programming methods, we compute the value of the game to be $25.9515$.
The Hider's optimal mixed strategy uses seven out of ten pure strategies---except for $(0,0,2,0)$, $(0,0,1,1)$, and $(0,0,0,2)$---with respective probabilities according to (\ref{eq:HiderEqual}), normalized to ensure it is a probability distribution.
\hfill $\Box$
\end{example}

\smallskip

The fact that the Hider's probabilities are chosen according to~(\ref{eq:HiderEqual}) is significant, and we conjecture that the optimal Hider strategy always has this property.
However, it is not clear which pure strategies to use in the optimal mixed strategy.
Based on the results in Section~\ref{sec:n=2} and Example~\ref{ex:CostEx1}, one may conjecture that optimal policy is to rank all Hider's pure strategies $\mathbf{x} =(x_1, \ldots, x_n)$ according to $\Pi_{i=1}^n c_i^{x_i}$, and then let the Hider's mixed strategy use $d$ pure strategies with largest such values, for some positive integer $d$.
We conclude this section with a counterexample to this conjecture, and leave the Searcher's optimal strategy as an open problem.

\begin{example}
\label{ex:CostEx2}
	Consider the multi-look search game with search cost, where $n=4$, $k=2$, and $(c_1, c_2, c_3, c_4) = (100, 10, 1, 0.99)$.
	Using linear programming methods, we compute the value of the game to be $201.0972$.
	The Hider's optimal mixed strategy uses the four pure strategies $(2,0,0,0)$, $(1,1,0,0)$, $(1,0,1,0)$, $(1,0,0,1)$, with respective probabilities according to Equation~(\ref{eq:HiderEqual}), normalized to ensure it is a probability distribution.
Please note that $c_2^2 =100 > 99 = c_1 c_4$, but the Hider's optimal strategy does not use $(0,2,0,0)$.
\hfill $\Box$
\end{example}


\section{Multi-look Search Game with Regret}
\label{sec:regret-multi}
This section presents a similar model to that in Section~\ref{sec:cost-multi}, but with a different objective function.
When the Searcher opens box $i$ and finds it empty, we say a \textit{regret} $c_i$ is incurred.
The total regret incurred by the Searcher is the sum of all the regret incurred in recovering all $k$ balls.
The Searcher wants to minimize the expected total regret, while the Hider wants to maximize it.


First, we prove that the Hider's mixed strategy in (\ref{eq:HiderEqual}) is still an equalizing strategy for this case, where the objective is the expected regret.

\begin{lemma}
For the muli-look search game with regret, against the Hider strategy $\mathbf{p}$ in (\ref{eq:HiderEqual}), any Searcher strategy has expected regret
\begin{equation}
V([n],k) \equiv T_1([n]) - \frac{T_{k+1}([n])}{T_k([n])}.
\label{eq:Vnk}
\end{equation}
\end{lemma}
\textit{Proof.}
This lemma follows from Lemma~\ref{le:Unk} and the fact that the expected regret $V([n],k)$ is given by
\begin{align*}
V([n],k) & = U([n],k) - \sum_{\mathbf x} p_{n,k}(\mathbf x) \sum_{j=1}^n x_j c_j \\
& = \frac{k T_{k+1}([n]) - \sum_{\mathbf x} \prod_{i=1}^n c_i^{x_i}  \sum_{j=1}^n x_j c_j}{T_k([n])},
\end{align*}
where the outer sum is over all Hider strategies $\mathbf x$. Hence, it is sufficient to show that
\[
(k+1)T_k([n]) - T_1([n]) T_k([n]) = \sum_{\mathbf x} \prod_{i=1}^n c_i^{x_i}  \sum_{j=1}^n x_j c_j.
\]
Without loss of generality, a typical term on the left-hand side can be written $c_1^{x_1} \cdots c_t ^{x_t}$, for some $1 \le t \le n$ with $x_1,\ldots,x_t >0$ and $x_1+\cdots + x_t = k+1$. Its coefficient is $(k+1)-t$.

We argue that its coefficient on the right-hand side is the same. The term appears $t$ times in the outer sum, with coefficients $x_1-1,\ldots,x_t-1$. For example, it appears with coefficient $x_1 - 1$ when the term $c_1^{x_1 - 1} c_2^{x_2}\cdots c_t^{x_t}(c_1(x_1-1)+c_2x_2+\cdots + c_tx_t)$ is expanded. The sum of these coefficients is $(x_1-1) + \cdots + (x_t - 1) = k+1 - t$, which completes the proof.
Interested readers can also prove this lemma by mathematical induction similar to the proof for Lemma~\ref{le:Unk}.
\hfill $\Box$



\smallskip

It turns out that the equalizing strategy in (\ref{eq:HiderEqual}) is the Hider's optimal strategy, and the value of the game is given in  (\ref{eq:Vnk}).
In Section~\ref{sec:Vnk}, we prove this result by presenting an optimal strategy for the Searcher that relies on a recursive algorithm.
In Section~\ref{sec:normal}, we present a few special cases, where the Searcher can achieve optimality with a simpler search strategy.

\subsection{An Optimal Search Strategy}
\label{sec:Vnk}

\begin{theorem}
Consider the muli-look search game with regret.
The value of the game with $n$ boxes and $k$ balls is given by $V([n], k)$ in (\ref{eq:Vnk}), for $n \geq 1$ and $k \geq 1$.
\label{th:Vnk}
\end{theorem}
\textit{Proof.}
We use mathematical induction on $n$ to prove this theorem.
For $n=1$, we have $V([1], k) = c_1 - c_1^{k+1} / c_1^k = 0$, which is clearly the value of the game.

Suppose the theorem is true when there are $n-1$ boxes, and consider the search game with $n$ boxes.
Let the Searcher deal with box $n$ first.
In particular, the Searcher decides the maximal number of times she wants to open box $n$, by choosing an integer $s \in \{0, 1, \ldots, k\}$, with the proviso that she will stop as soon as a search in box $n$ finds no ball.
Let $t$ denote the number of balls hidden in box $n$.
If the Searcher recovers $t < s$ balls, then a cost $c_n$ is incurred on the $(t+1)$th search in box $n$, and the Searcher knows there are $k-t$ balls hidden among boxes $1,\ldots, n-1$, so the expected total regret becomes $c_n + V([n-1], k-t)$, according to the induction hypothesis.
If the Searcher recovers $s$ balls, then she puts box $n$ aside, and uses the  policy that would be optimal on boxes $1, 2, \ldots, n-1$ assuming $n-s$ balls are hidden among those $n-1$ boxes.
The Searcher will not open box $n$ again, unless it turns out that fewer than $n-s$ balls are hidden among boxes $1,\ldots, n-1$, after incurring a total regret $\sum_{i=1}^{n-1} c_i$.
If $t=s$, then the expected total regret becomes $V([n-1], k-s)$.
If $t>s$, then the total regret becomes $\sum_{i=1}^{n-1} c_i$.
The payoff matrix is given in Table~\ref{tab:adaptive}, where the Hider chooses row $t$ (the number of balls to hide in box $n$), and the Searcher chooses column $s$ (the maximal number of times to open box $n$), for $s,t = 0, 1, \ldots, k$.

\begin{table}[h!]
	\caption{The matrix game, where the Hider choose the number of balls to hide in box $n$ (row), and the Searcher decides the maximal number of times to open box $n$ (column).}
	\label{tab:adaptive}
	\begin{center}
		\begin{scriptsize}
			\begin{tabular}{c|cccccc}
				& $0$ & $1$ & $\cdots$ & $k-2$ & $k-1$ & $k$ \\
				\hline
				0 & $V([n-1], k)$ & $c_n + V([n-1], k)$ & $\cdots$  & $c_n + V([n-1], k)$  & $c_n + V([n-1], k)$  & $c_n + V([n-1], k)$\\
				1 & $\sum_{i=1}^{n-1} c_i$ & $V([n-1], k-1)$ & $\cdots$  & $c_n + V([n-1], k-1)$ & $c_n + V([n-1], k-1)$ & $c_n + V([n-1], k-1)$ \\
				$\vdots$ & $\cdots$ & $\cdots$ & $\cdots$  &  \\
				$k-2$ & $\sum_{i=1}^{n-1} c_i$ & $\sum_{i=1}^{n-1} c_i$ & $\cdots$  & $V([n-1], 2)$ & $c_n + V([n-1], 2)$ & $c_n + V([n-1], 2)$\\
				$k-1$ & $\sum_{i=1}^{n-1} c_i$ & $\sum_{i=1}^{n-1} c_i$ & $\cdots$ & $\sum_{i=1}^{n-1} c_i$ & $V([n-1], 1)$ & $c_n + V([n-1], 1)$ \\
				$k$ & $\sum_{i=1}^{n-1} c_i$ & $\sum_{i=1}^{n-1} c_i$ & $\cdots$  & $\sum_{i=1}^{n-1} c_i$ & $\sum_{i=1}^{n-1} c_i$ & $0$
			\end{tabular}
		\end{scriptsize}
	\end{center}
\end{table}

Now, we claim that the value of the matrix game in Table~\ref{tab:adaptive} is given by $V([n],k)$ in (\ref{eq:Vnk}), and the optimal mixed strategy for the Searcher is to use column $s$ with probability
\begin{equation}
P_{s} = \frac{\frac{T_{k-s}([n])}{T_{k-s+1}([n-1])} - \frac{T_{k-s-1}([n])}{T_{k-s}([n-1])}}{\frac{T_k([n])}{T_{k+1}([n-1])}},
\label{eq:MultiRegretOpt}
\end{equation}
for $s=0, 1, \ldots, k$, where we define
\[
T_{-1}([n]) \equiv 0, \qquad T_0([n]) \equiv 1,
\]
for notational convenience.
Recall that the Hider's mixed strategy in (\ref{eq:HiderEqual}) yields expected regret $V([n], k)$ for any search strategy.
To prove this claim, it remains to show that the Searcher's mixed strategy in  (\ref{eq:MultiRegretOpt}) yields the same expected regret $V([n],k)$ against row $t$, for $t=0,\ldots,k$.

Suppose that the Hider hides $t$ balls in box $n$.
The Searcher's expected regret with the mixed strategy given in  (\ref{eq:MultiRegretOpt}) is equal to
\begin{align*}
&\left(\sum_{i=1}^{n-1} c_i \right) \left( \sum_{i=0}^{t-1} P_i \right) + V([n-1], k-t) \cdot P_{t} + (c_n + V([n-1], k-t)) \left( \sum_{i=t+1}^{k} P_i \right) \\
&= \left(\sum_{i=1}^{n-1} c_i \right) \left( \sum_{i=0}^{t-1} P_i \right) + V([n-1], k-t) \left( \sum_{i=t}^{k} P_i \right) + c_n \left( \sum_{i=t+1}^{k} P_i \right).
\end{align*}
Using the induction hypothesis $V([n-1],k-t) = T_1([n-1]) - \frac{T_{k-t+1}([n-1])}{T_{k-t}([n-1])}$, and then using  (\ref{eq:MultiRegretOpt}), the preceding becomes
\begin{align*}
&\sum_{i=1}^{n-1} c_i  -  \frac{T_{k-t+1}([n-1])}{T_{k-t}([n-1])}  \left( \sum_{i=t}^{k} P_i \right) + c_n  \left( \sum_{i=t+1}^{k} P_i \right) \\
&= \sum_{i=1}^{n} c_i - c_n -  \frac{T_{k-t+1}([n-1])}{T_{k-t}([n-1])}  \left( \frac{\frac{T_{k-t}([n])}{T_{k-t+1}([n-1])}}{\frac{T_{k}([n])}{T_{k+1}([n-1])} } \right) + c_n  \left( \frac{\frac{T_{k-t-1}([n])}{T_{k-t}([n-1])}}{\frac{T_{k}([n])}{T_{k+1}([n-1])} } \right) \\
&= T_1([n]) - \frac{T_{k+1}([n-1])}{T_{k}([n])} \left( \frac{c_n \cdot T_{k}([n])}{T_{k+1}([n-1])} + \frac{T_{k-t}([n])}{T_{k-t}([n-1])} -  \frac{c_n \cdot T_{k-t-1}([n])}{T_{k-t}([n-1])} \right).
\end{align*}
Using the identity
\[
c_n \cdot T_{k-t-1}([n]) + T_{k-t}([n-1]) = T_{k-t}([n]),
\]
we end up with
\begin{align*}
&T_1([n]) - \frac{T_{k+1}([n-1])}{T_{k}([n])} \left(  \frac{c_n \cdot T_{k}([n])}{T_{k+1}([n-1])} + 1 \right) \\
&= T_1([n]) - \frac{c_n \cdot T_{k}([n])+ T_{k+1}([n-1])}{T_{k}([n])}  \\
&= T_1([n]) - \frac{T_{k+1}([n])}{T_{k}([n])},
\end{align*}
which proves the value of the matrix game in Table~\ref{tab:adaptive} is $V([n], k)$.
Consequently, the theorem holds true for $n$ boxes, which completes the proof.
\hfill $\Box$

\smallskip

A byproduct of the proof of Theorem~\ref{th:Vnk} is an optimal search strategy, which relies on a recursive algorithm.

\smallskip

\noindent
\textbf{Search}$([n], k)$
\begin{enumerate}
\item
Generate a random number, denoted by $s$, according to the probability mass function in (\ref{eq:MultiRegretOpt}).
\item
Open box $n$ repeatedly until one of the following:
\begin{enumerate}
\item
Recover $s$ balls: Run \textbf{Search}$([n-1], k-s)$ on boxes $1, 2, \ldots, n-1$.
If the Searcher finds fewer than $k-s$ balls in boxes $1, 2, \ldots, n-1$ (after incurring regret $\sum_{i=1}^{n-1} c_i$), then open box $n$ again, repeatedly, until recovering all $k$ balls.
\item
Recover $t < s$ balls and incur regret $c_n$:
Run \textbf{Search}$([n-1], k-t)$ on boxes $1, 2, \ldots, n-1$.
\end{enumerate}
	
\end{enumerate}

\subsection{Normal Search Strategies in Special Cases}
\label{sec:normal}
In this subsection, we discuss a more restricted class of strategies for the Searcher, which we call {\em normal strategies}. Roughly speaking, a normal strategy sets out a search plan in advance and does not allow the Searcher to deviate from it by using information acquired during the search. The advantage of normal strategies is that they can be described concisely, and only rely on the search history at every stage in a limited way.

A \textit{search sequence} is a predetermined sequence of box numbers, such that the Searcher opens boxes following the sequence, with the proviso that she will skip a box if she already found out that the box is empty.
Specifically, it is a sequence of box numbers with length $n \cdot k$, which contains each box number $i=1,\ldots,n$ precisely $k$ times. The total number of different search sequences is $(nk)!/(k!)^n$. A \textit{normal strategy} is a mixed strategy of these search sequences. For a search sequence $a$ and a Hider strategy $\mathbf{x}$, we write $R(a,\mathbf{x})$ for total regret incurred when the Searcher follows $a$ and the Hider follows $\mathbf{x}$.

For example, suppose $n=3$ and $k=2$, and consider the search sequence $a=\{1,3,3,2,1,2\}$. For the Hider strategy $\mathbf{x}=(0,1,1)$, the Searcher will first open box $1$, then she will open box $3$ twice, and finally she will open box $2$, giving a total regret of $R(a,\mathbf{x})=c_1+c_3$.
For the Hider strategy $\mathbf{x}=(2,0,0)$, the Searcher will open box $1$, then box $3$, then box $2$ and finally box $1$ again, with a regret of $R(a,\mathbf{x})=c_2+c_3$.

We make the observation that the order of the first $k$ terms in a search sequence is irrelevant, since the balls cannot be found until at least $k$ searches have been performed.
We also show that the order of the last $k$ terms in a search sequence is irrelevant.

\begin{lemma} \label{lem:normal}
For the muli-look search game with regret, changing the order of the first or last $k$ terms in a search sequence does not change the payoff of this search against any Hider strategy.
\end{lemma}
\textit{Proof.}
We have already noted that the order of the first $k$ terms is irrelevant.
To show that the order of the last $k$ terms is irrelevant, write $y_i$ for the number of searches in box $i$ in the last $k$ terms, where $y_1+\cdots + y_n=k$.
Consider some fixed Hider strategy $(x_1,\ldots,x_n)$.
	
First, suppose that $x_i +y_i \le k$ for all $i$.
It is straightforward to see that all the balls will be found before the last $k$ terms of the search, so the order of the last $k$ terms is irrelevant.
	
Second, suppose that there is some $i$ for which $x_i+y_i > k$.
For $j \neq i$, we have that
\begin{align}
x_j + y_j &= \left(k - \sum_{t \neq j} x_t \right) + \left(k - \sum_{t \neq j} y_t \right) \nonumber \\
	&\le (k - x_i) + (k-y_i) \nonumber \\
	& < k. \label{eq:xy}
\end{align}
At some point before the last $k$ terms of the search, the Searcher has looked for a ball in box $j$ and not found it, so she knew box $j$ has become empty, for $j \neq i$.
Hence, she will look only in box $i$ in the last $k$ terms, so the order of the last $k$ terms is irrelevant.
\hfill $\Box$

\smallskip

We say that the search sequence starts (or ends, respectively) with $(y_1,\ldots,y_n)$, with $\sum_{i=1}^n y_i = k$, where $y_i$ denotes the times box $i$ shows up in the first (or last, respectively) $k$ terms in the search sequence.
We make another simple observation.

\begin{lemma} \label{lem:symm2}
Consider the muli-look search game with regret.
Let $\mathbf{x} = (x_1, \ldots, x_n)$ and $\mathbf{y} = (y_1, \ldots, y_n)$ denote two  Hider strategies, and let $a_{\mathbf{x}}$ and $a_{\mathbf{y}}$ be search sequences that end with $\mathbf{x}$ and $\mathbf{y}$, respectively. 
If there exists some $i$ such that $x_i+y_i > k$, then the regret function is symmetric in $\mathbf{x}$ and $\mathbf{y}$.
That is, $R(a_{\mathbf{y}},\mathbf{x})=R(a_\mathbf{x},\mathbf{y})$.
\end{lemma}
\textit{Proof.}
It follows from (\ref{eq:xy}) that $x_j + y_j < k$ for $j \neq i$.
It is straightforward to verify that $R(a_{\mathbf{y}},\mathbf{x})=R(a_\mathbf{x},\mathbf{y}) = \sum_{j \neq i} c_j$.
\hfill $\Box$

\smallskip

In order to prove the optimality of certain normal search strategies in this section, we will use a simple result also used in Lidbetter (2013), which we repeat here without proof.

\begin{lemma} \label{lem:symm} 
Consider a two-player, zero-sum game between Players I and II with payoff matrix $M$, corresponding to Player I's payoffs. If $M$ is symmetric and $\mathbf{x}^*$ is a mixed strategy for Player I that
makes Player II indifferent between all his strategies, then the strategy pair $(\mathbf{x}^*, \mathbf{x}^*)$ forms an equilibrium.
\end{lemma}

We now consider three special cases, and present an optimal \textit{normal} Searcher strategy in each case.

\subsubsection{Special Case: $n=2$}
\label{sec:k=1}
When there are only $n=2$ boxes, the Hider has $k+1$ pure strategies $(x_1, x_2)$, where $x_1+x_2=k$, for $x_1=0, 1, \ldots, k$.
Due to Lemma~\ref{lem:normal}, a search sequence can be delineated by $(y_1,y_2)$, where $y_i$ is the number of times that box $i$ appears in the last $k$ terms of the sequence, with $i=1,2$ and $y_1+y_2 = k$.
In the notation of Lemma~\ref{lem:symm2}, the strategy $a_{\mathbf{y}}$ is uniquely specified, and it is straightforward to verify that
\[
R(a_{\mathbf{y}},\mathbf{x}) = R(a_\mathbf{x},\mathbf{y}) =
\left\{
\begin{array}{ll}
0, & \mbox{if } x_1+y_1=k,\\
c_1, & \mbox{if } x_1+y_1<k,\\
c_2, & \mbox{if } x_1+y_1>k.
\end{array}
\right.
\]
Since the payoff matrix is symmetric and the Hider strategy $\mathbf{p}$ in (\ref{eq:HiderEqual}) makes the Searcher indifferent between all her strategies, the same mixed strategy $\mathbf{p}$ is optimal for both players, by Lemma~\ref{lem:symm}.
That is, the Hider and Searcher should both choose a pair $(y_1,y_2)$ with probability $c_1^{y_1} c_2^{y_2} / T_k([2])$, which is consistent with the probability in (\ref{eq:MultiRegretOpt}).

\subsubsection{Special Case: $k=1$}
\label{sec:k=1}
Consider the case with $k=1$ ball.
The Hider has $n$ pure strategies, corresponding to where to he hides the ball.
There are $n!$ different search sequences, each corresponding to a permutation of the $n$ boxes.
Write $\pi_j$ for the Searcher's normal strategy that chooses each search sequence that ends with box $j$ with equal probability $1 / (n-1)!$, for $j = 1, \ldots, n$.
Write $R(\pi_j, i)$ for the expected regret if the Searcher uses the normal strategy $\pi_j$, while the Hider hides the ball in box $i$, for $i, j = 1, \ldots, n$.
One can check that
\[
R(\pi_{j}, i) = R(\pi_{i}, j) =
\left\{
\begin{array}{ll}
\sum_{t \neq i } c_t, & \mbox{if } i=j,\\
\frac{1}{2} \sum_{t \neq i,j} c_t, & \mbox{if } i \neq j.
\end{array}
\right.
\]
Hence, the payoff matrix is symmetric between the two players.
Since the Hider strategy of equation (\ref{eq:HiderEqual}) makes the Searcher indifferent, by Lemma~\ref{lem:symm}, we conclude that the Searcher's optimal strategy is to use $\pi_j$ with probability $c_j / \sum_{t=1}^n c_t$.

\subsubsection{Special Case: $k=2$}
Consider the special case with $k=2$ balls.
For every $\mathbf{y}$ with $y_1+\cdots+y_n=2$, we will define a normal strategy $T_\mathbf{y}$ that ends with $\mathbf{y}$, such that the expected regret $R(T_\mathbf{y},\mathbf{x})$ is symmetric in $\mathbf{x}$ and $\mathbf{y}$.
To describe the search subsequence before the last two searches, we write $\sigma(A)$ for a uniformly random permutation of all elements in a set $A$.
For example, $\sigma(3, 4, \ldots,n)\sigma(1, 2, \ldots,n)$ means a random permutation of boxes $\{3,4,\ldots,n\}$ followed by an independent random permutation of boxes $\{1, 2,\ldots,n\}$.
For another example, $\sigma(22, 33, \ldots,nn)$ means that, if the Searcher opens a box and finds a ball, then she will immediately open the same box again, while she follows a random order of boxes $\{2, 3, \ldots,n\}$ in doing so.



With $k=2$ balls, a search sequence has length $2n$ with each box appearing twice in the sequence, and a normal strategy is a mixed strategy on these search sequences.
We next give a normal strategy that ensures the payoff function is symmetric.
To describe a normal strategy, we need to specify what the Searcher does before her last two searches, which is either $\mathbf{y}=\{2,0,\ldots,0\}$, or $\mathbf{y}=\{1,1,0,\ldots,0\}$, or some permutation of these vectors.

\begin{definition} [A normal search strategy for $k=2$]
\label{de:normal}
If the Searcher ends with $\mathbf{y} = (2,0,\ldots,0)$, then she starts with $\sigma(22, 33, \ldots,nn)$.
If the Searcher ends with $\mathbf{y}= (1,1,0,\ldots,0)$, then she starts with each of the following with probability $1/3$:
\[
\sigma(3,4,\ldots,n),\sigma(1,2,\ldots,n); \qquad 1,\sigma(33,44,\ldots,nn),2; \qquad 2,\sigma(33,44,\ldots,nn),1.
\]
\end{definition}

\begin{lemma}
For the muli-look search game with regret, write $\mathbf{x}$ for a Hider's strategy, and $T_\mathbf{y}$ for the Searcher's normal strategy in Definition~\ref{de:normal}.
The payoff function $R(T_\mathbf{y},\mathbf{x})$ is symmetric in $\mathbf{x}$ and $\mathbf{y}$; that is, $R(T_\mathbf{y},\mathbf{x}) = R(T_\mathbf{x},\mathbf{y})$.
\label{le:Ty}
\end{lemma}
\textit{Proof.}
Up to symmetry, there are 7 cases:
\begin{enumerate}
\item
$\mathbf{y} =\mathbf{x} = (2,0,\ldots,0)$.
By definition, $R(T_\mathbf{y},\mathbf{x}) = R(T_\mathbf{x},\mathbf{y})$.

\item
$\mathbf{y} =(2,0,\ldots,0)$ and $\mathbf{x} = (0,2,0,\ldots,0)$.
The expected regret is
\[
R(T_\mathbf{y},\mathbf{x}) = \frac{1}{2} (c_3 + c_4 + \cdots + c_n) = R(T_\mathbf{x},\mathbf{y}).
\]

\item
$\mathbf{y} =(2,0,0,\ldots,0)$ and $ \mathbf{x} = (1,1,0,\ldots,0)$.
By Lemma~\ref{lem:symm2}, $x_1 + y_1 = 3 > 2 = k$, so $R(T_\mathbf{y},\mathbf{x}) = R(T_\mathbf{x},\mathbf{y})$.

\item
$\mathbf{y} =(2,0,0,\ldots,0)$ and $ \mathbf{x} = (0,1,1,0,\ldots,0)$. The expected regret is
		\[
		R(T_\mathbf{y},\mathbf{x})=\frac{1}{2}(c_2+c_3)+\frac{2}{3} (c_4+\cdots+c_n).
		\]
		Swapping the strategies, compute
		\begin{align*}
		R(T_\mathbf{x},\mathbf{y})&= \frac{1}{3} \left( (c_4+\cdots+c_n)+ \frac{1}{2} (c_2+c_3) \right) + \frac{2}{3} \left( \frac{1}{2} (c_2+c_3+c_4+\cdots+c_n) \right)\\
		&= \frac{1}{2} (c_2+c_3)+ \frac{2}{3}(c_4+\cdots+c_n).
		\end{align*}

\item
$\mathbf{y}= \mathbf{x}=(1,1,0,\ldots,0)$. 
By definition, $R(T_\mathbf{y},\mathbf{x}) = R(T_\mathbf{x},\mathbf{y})$.

\item
$\mathbf{y}=(1,1,0,\ldots,0)$ and $\mathbf{x}=(0,1,1,0,\ldots,0)$. The expected regret is
		\begin{align*}
		R(T_\mathbf{y},\mathbf{x})&= \frac{1}{3} \left( \frac{1}{2} (c_1+c_3) + (c_4+\cdots+c_n ) \right) + \frac{1}{3} (c_1+c_3+c_4+\cdots+c_n) \\
		& \quad + \frac{1}{3} \left( \frac{1}{2} (c_4+\cdots+c_n) \right)\\
		&= \frac{1}{2} (c_1+c_3)+ \frac{5}{6} (c_4+\cdots+c_n).
		\end{align*}
		By symmetry, $R(T_\mathbf{y},\mathbf{x})=R(T_\mathbf{x},\mathbf{y})$.
\item 
$\mathbf{y} =(1,1,0,\ldots,0)$ and $\mathbf{x} =(0,0,1,1,0,\ldots,0)$.
The expected regret is
		\begin{align*}
		R(T_\mathbf{y},\mathbf{x}) &= \frac{1}{3} \left( \frac{2}{3} (c_5+\cdots+c_n)\right) + \frac{2}{3} \left(\frac{1}{2}(c_1+c_2)+ \frac{1}{2} (c_3+c_4)+ \frac{2}{3}(c_5+\cdots+c_n) \right)\\
		&= \frac{1}{3} (c_1+c_2+c_3+c_4) +  \frac{2}{3} (c_5+\cdots+c_n),
		\end{align*}
		By symmetry, $R(T_\mathbf{y},\mathbf{x})=R(T_\mathbf{x},\mathbf{y})$.
\end{enumerate}
Since $R(T_\mathbf{y},\mathbf{x})=R(T_\mathbf{x},\mathbf{y})$ in all cases, the proof is completed.
\hfill $\Box$

Finally, it follows from Lemma~\ref{lem:symm} and Lemma~\ref{le:Ty} that the normal strategy in Definition~\ref{de:normal} is optimal for the Searcher for the multi-look search game with regret, in the special case where there are $k=2$ balls.

We leave as an open problem the question of whether the Searcher always has an optimal normal strategy in this game.

\section{Single-look Search Game with Regret}
\label{sec:regret-single}
This section presents a similar search problem to that in Section~\ref{sec:regret-multi}, with the exception that the Hider can hide at most one ball in a box.
Therefore, the Hider chooses $k$ boxes to hide the $k$ balls, and the Searcher chooses a search sequence, which is a permutation of the $n$ boxes.
The Hider has ${n \choose k}$ pure strategies, while the Searcher has $n!$ pure strategies.
The Searcher wants to minimize the expected total regret, while the Hider wants to maximize it.
The problem can also be viewed as a variation to that in \cite{lidbetter}, but the objective here is the total regret rather than the total search cost.

Let $S_k([n])$ denote the sum over all products $\Pi_{i \in A} c_i$, for all $A$ that is a subset of $\{1,\ldots, n\}$ of size $k$.
In other words,
\[
S_k([n]) = \sum_{A \subseteq \{1,\ldots, n\}, |A|=k} \left( \prod_{i \in A} c_i \right).
\]
For instance, $S_2([3]) = c_1 c_2 + c_2 c_3 + c_1 c_3$.
Write a Hider's pure strategy as $\mathbf{x} = (x_1, x_2, \ldots, x_n)$, where $x_i=1$ if a ball is hidden in box $i$, or $x_i=0$ otherwise, with $\sum_{i=1}^n x_i = k$.
Consider the Hider strategy, which chooses a pure strategy $\mathbf{x} = (x_1, \ldots, x_n)$ with probability
\begin{equation}
p_{n, k}(\mathbf{x}) =  \frac{\prod_{i=1}^n c_i^{x_i}}{S_k([n])}.
\label{eq:SingleHiderEqual}
\end{equation}
Our first result is that the Hider's mixed strategy in (\ref{eq:SingleHiderEqual}) is an equalizing strategy. The proof is similar to that of Lemma~\ref{le:Unk}, and can be found in the Appendix.

\begin{lemma}
\label{le:Wnk}
For the single-look search game with regret, against the Hider strategy $\mathbf{p}$ in (\ref{eq:SingleHiderEqual}), any Searcher strategy has expected regret
\begin{equation}
W([n],k) \equiv \frac{k S_{k+1}([n])}{S_k([n])}.
\label{eq:Wnk}
\end{equation}
\end{lemma}




The Hider's equalizing strategy in (\ref{eq:SingleHiderEqual}) is the optimal strategy in a few special cases.
For example, when $k=1$, the model reduces to the one in Section~\ref{sec:k=1}, since it makes no difference how many balls can be hidden in one box if there is only 1 ball.
For another example, when $c_i=c$ for $i=1,\ldots, n$, then due to symmetry, it is optimal for both the Hider and the Searcher to respectively choose each pure strategy uniformly random.
In general, however, the equalizing strategy in (\ref{eq:SingleHiderEqual}) is not necessarily an optimal strategy.
We present the solution to a special case $k = n-1$ in Section~\ref{sec:k=n-1}, and then offer a few observations in Section~\ref{sec:single-equal-general}.

\subsection{Special Case: $k=n-1$}
\label{sec:k=n-1}
Suppose there are $n$ boxes and $k=n-1$ balls.
The Hider has $n$ pure strategies, by deciding which one box to leave empty.
The Searcher essentially needs to guess which box is empty.
Hence, the Searcher also has $n$ pure strategies, by deciding which box to search last.
If the Hider leaves box $i$ empty, and the Searcher searches box $j$ last, then the total regret is
\[
R(i,j) = \left\{
\begin{array}{ll}
0, & \mbox{ if } i = j, \\
c_i, & \mbox{ if } i \neq j.
\end{array}
\right.
\]

%

\begin{lemma}
\label{le:RegretEqualizing}
Consider the single-look search game with regret.
In the case $k=n-1$, if
\begin{equation}
q_{n}(j) = 1 - \frac{(n-1) /c_j}{\sum_{i=1}^n 1/c_i} \geq 0,
\label{eq:qnj}
\end{equation}
for $j=0, \ldots, n$, then the Searcher's optimal mixed strategy is to search box $j$ last with probability $q_n(j)$.
The Hider's optimal mixed strategy is given in (\ref{eq:SingleHiderEqual}), and the value of the game is
\begin{equation}
\frac{(n-1) S_{n}([n])}{S_{n-1}([n])}.
\label{eq:k=n-1,value}
\end{equation}
\end{lemma}
\textit{Proof.}
First, it is straightforward to check $\sum_{j=1}^n q_n(j) = 1$.
Next, if the Hider leaves box $i$ empty, then the Searcher's mixed strategy in (\ref{eq:qnj}) yields expected regret
\begin{align*}
\sum_{j=1}^n q_n(j) R(i,j) &= c_i (1- q_n(i) ) 
= c_i \times \frac{(n-1) /c_i}{\sum_{i=1}^n 1/c_i}
= \frac{(n-1) S_{n}([n])}{S_{n-1}([n])},
\end{align*}
where the last equality follows by multiplying $S_n([n]) = \Pi_{i=1}^n c_i$ to the numerator and the denominator.
Since the Hider and the Searcher each has a mixed strategy that guarantees the same payoff in (\ref{eq:k=n-1,value}) regardless of what the other does, the proof is completed.
\hfill $\Box$

\smallskip

If $q_n(j) < 0$ for some $j$, then the optimal strategy is more complicated.
Without loss of generality, for the remainder of this section we assume that $c_1 \geq c_2 \geq \cdots \geq c_n$. 
Since $q_n(j)$ decreasaes in $j$, it is a legitimate probability distribution if $q_n(n) \geq 0$, or equivalently,
\[
1 - \frac{(n-1) /c_n}{\sum_{i=1}^n 1/c_i} \geq 0,
\]
or equivalently,
\begin{equation}
c_n \geq \frac{n-2}{\sum_{i=1}^{n-1} 1 / c_i}.
\label{eq:criterion}
\end{equation}
Otherwise, if (\ref{eq:criterion}) does not hold, then $q_n(n) < 0$, which implies that the Searcher should never save box $n$ to check last.
Intuitively, if $c_n$ is too small, then it is not worthwhile for the Hider to leave box $n$ empty so as to collect a potential small regret.
Hence, it is optimal for the Hider to always put a ball in box $n$, and for the Searcher to always open box $n$ first.
For example, suppose $c_1 = c_2 = 100$ and $c_3=1$.
The Searcher's optimal strategy is to open box 3 first, and then open boxes 1 or 2 with equal probability.
The Hider's optimal strategy is to put 1 ball in box 3, and the other ball either in box 1 or box 2 with probability 0.5.
The value of the game is the expected regret $0.5 \times 100 = 50$.
If the Hider does not put a ball in box 3, but 1 ball each in boxes 1 and 2 instead, then the expected regret is only 1, when the Searcher searches box 3 first.

Intuitively, we can use equation (\ref{eq:criterion}) to determine how the Hider divides the $n$ boxes into two groups.
If the Hider chooses some $m \geq 2$ and puts 1 ball each in boxes $m+1, m+2, \ldots, n$, and play a game with $m-1$ balls and boxes $1, 2, \ldots, m$, then the Searcher's equalizing strategy guarantees expected regret $W([m], m-1)$.
A sensible strategy is for the Hider to choose the best $m$ to maximize $W([m], m-1)$, so the value of the game would be
\[
\max_{m=2, \ldots, n} W([m], m-1) = \max_{m=2, \ldots, n} \frac{(m-1) S_m([m])}{S_{m-1}([m])}.
\]
We present a lemma before proving this result.

\begin{lemma}
For $c_1 \geq c_2 \geq \ldots \geq c_n$, define
\begin{equation}
a_m = \sum_{i=1}^{m-1} \frac{1}{c_i} - \frac{m-2}{c_m},
\label{eq:a_m}
\end{equation}
for $m=2, \ldots, n$.
The sequence $a_2, a_3, \ldots, a_n$ is decreasing (weakly).
\end{lemma}
\textit{Proof.}
Take the difference
\[
a_{m+1} - a_m = \frac{1}{c_m} - \frac{m-1}{c_{m+1}} + \frac{m-2}{c_m} = (m-1) \left( \frac{1}{c_m} - \frac{1}{c_{m+1}} \right) \leq 0,
\]
since $c_{m+1} \leq c_m$, which completes the proof.
\hfill $\Box$

\smallskip

\begin{theorem}
Consider the single-look search game with regret, and suppose $k=n-1$. 
Define $a_m$ as in equation (\ref{eq:a_m}), for $m=2, \ldots, n$, and let $b = \max \{m: a_m \geq 0\}$.
The Hider's optimal strategy is to first put 1 ball each in boxes $b+1, b+2, \ldots, n$, and then hide the other $b-1$ balls among boxes $1, 2, \ldots, b$ according to (\ref{eq:SingleHiderEqual}).
The Searcher's optimal strategy is to first open boxes $b+1, b+2, \ldots, n$.
If any of these boxes is empty, then open boxes $1, 2, \ldots, b$ to find all balls; otherwise, use the mixed strategy according to $q_b(j)$ in (\ref{eq:qnj}) for boxes $1, 2, \ldots, b$.
The value of the game is $W([b], b-1)$.
\end{theorem}
\textit{Proof.}
With the stated Hider strategy, the Hider guarantees an expected regret of $W([b], b-1)$ regardless of what the Searcher does (Lemma~\ref{le:Wnk}), so it remains to show that the stated Searcher strategy guarantees an expected regret at most $W([b], b-1)$ for any Hider strategy.

Suppose the Searcher uses the stated strategy, and opens boxes $b+1, \ldots, n$ first.
If any of these boxes is empty, then the regret is at most $c_{b+1}$.
Otherwise, the expected regret is $W([b], b-1)$.
Since $a_{b+1} < 0$, it follows that
\[
\sum_{i=1}^{b} \frac{1}{c_i} < \frac{b-1}{c_{b+1}},
\]
or equivalently,
\[
c_{b+1} <  \frac{(b-1)}{\sum_{i=1}^{b} c_i} = W([b], b-1),
\]
which completes the proof.
\hfill $\Box$

\subsection{The General Case}
\label{sec:single-equal-general}

Because each box can contain just one ball, putting a ball in a box takes away the opportunity for the Hider to collect a regret from that box, while leaving a box empty allows such an opportunity.
To maximize the total regret, the Hider needs to \textit{allocate} these $n-k$ opportunities carefully.
If $c_i$ is too small, then it is intuitive that the Hider should not waste an opportunity there to collect a very small regret.
The optimal policy in Section~\ref{sec:k=n-1} leads to a conjecture of the optimal policy.
The Hider puts one ball each in the $d$ boxes having the smallest search costs, and distributes the remaining $k-d$ balls among the other $n-d$ boxes according to the equalizing strategy in (\ref{eq:SingleHiderEqual}).
This conjecture, however, is not true, as seen in the following example.

\begin{example}
\label{ex:RegretEx1}
Consider the single-look search game with regret, where $n=4$, $k=2$, and $(c_1, c_2, c_3, c_4) = (100, 10, 1, 0.99)$.
Using linear programming methods, we compute the value of the game to be $10.0405$.
The Hider's optimal mixed strategy uses five out of six pure strategies---except for $(1,1,0,0)$---with respective probabilities according to (\ref{eq:SingleHiderEqual}), normalized to ensure it is a probability distribution.
\hfill $\Box$
\end{example}

\smallskip

In the single-look search problem with regret, we again find that the Hider's optimal mixed strategy always chooses probabilities according to (\ref{eq:SingleHiderEqual}).
However, it is not clear which pure strategies to use in the optimal mixed strategy.
One conjecture is to rank all Hider's pure strategies $\mathbf{x} =(x_1, \ldots, x_n)$ according to $\Pi_{i=1}^n c_i^{x_i}$, and then let the Hider's mixed strategy use $d$ pure strategies with smallest such values, for some positive integer $d$.
Although this conjecture is consistent with the results in Section~\ref{sec:k=n-1} and Example~\ref{ex:RegretEx1}, we give a counterexample below.

\begin{example}
\label{eq:RegretEx2}
Consider the single-look search game with regret, where $n=4$, $k=2$, and $(c_1, c_2, c_3, c_4) = (100, 10, 9.9, 1)$.
Using linear programming methods, we compute the value of the game to be $17.4229$.
The Hider's optimal mixed strategy uses $(1,0,0,1)$, $(0,1,0,1)$, $(0,0,1,1)$, with respective probabilities according to (\ref{eq:SingleHiderEqual}), normalized to ensure it is a probability distribution.
Please note that $c_2 c_3 =99 < 100 = c_1 c_4$, but the Hider's optimal strategy does not use $(0,1,1,0)$.
\hfill $\Box$
\end{example}

It appears rather difficult to derive the Hider's optimal policy in general.
In all of our numerical examples, we observe that the Hider's optimal policy chooses a subset of her pure strategies, and uses each pure strategy in the subset with probability proportional to that in (\ref{eq:SingleHiderEqual}).
It is plausible that the Hider's optimal strategy always obeys this rule, since it does not allow the Searcher to learn anything about the Hider's strategy in the search process.
We leave this observation as a conjecture.

\section{Conclusion}
\label{sec:conclusion}
This paper concerns a search game in which the Hider distributes $k$ balls among $n$ boxes, and the Searcher needs to find all balls by opening the boxes one by one.
There are four versions of this search problem, depending on whether the Searcher finds just one ball or all balls in the box when it is opened, and whether the objective is the expected total search cost or the expected total regret.
As summarized in Table~\ref{tab:results}, our work complements that in \cite{lidbetter}, and contributes to three versions of this search problem.

Remarkably, the Hider's equalizing strategy---in (\ref{eq:HiderEqual}) for the multi-look game and in (\ref{eq:SingleHiderEqual}) for the single-look game---does not depend on whether the objective function is the expected total cost or the expected total regret.
In addition, the Hider's equalizing strategy is the optimal strategy in many cases.
In all examples we examine, where the Hider's equalizing strategy is not optimal, the optimal Hider's optimal strategy is to use a subset of her pure strategies with probabilities proportional to those in her equalizing strategy. 
Call this property of a strategy the {\em equalizing property}.
We conjecture that the Hider's optimal strategy always has the equalizing property, since in this case, whenever the Hider opens a box, the remaining hidden balls will be hidden according to a strategy with the equalizing property.

Besides the open problems posed in Sections~\ref{sec:cost-multi}, \ref{sec:regret-multi} and \ref{sec:regret-single}, there are a few future research directions.
This paper assumes that the Searcher will always find a ball in a nonempty box, but in practice, the Searcher may fail to find a ball in the box she opens.
If the search has a deadline, then the Searcher's objective becomes finding as many objects as possible.
In some applications, there may exist geographical constraints on where the Searcher can go next from her current location.
All these extensions are important in practice, but require substantial new mathematical formulation.

\bibliographystyle{apalike}
\bibliography{references}

\section*{Appendix}

\textit{Proof of Lemma~\ref{le:Wnk}.}
The proof is by induction on $n+k$. For $n=k=1$, we have 
\[
\frac{1 \cdot S_{2}([1])}{S_1([1])} =  \frac{1 \cdot 0}{c_1}= 0,
\]
which is clearly the expected regret, since the Searcher will open box 1 exactly once to find the only ball, and be done with it.

Suppose $n+k  \geq 3$ and suppose without loss of generality that the Searcher starts by opening box $n$.
As noted already, whether or not a ball is found, the remaining balls will be hidden according to $\mathbf{p}$ in (\ref{eq:SingleHiderEqual}).
More precisely, if a ball is found (which happens with probability $c_n S_{k-1}([n-1])/S_k([n])$, the remaining balls are hidden according to $p_{n,k-1}$; if a ball is not found (which happens with probability $S_k([n-1])/S_k([n])$, the remaining balls are hidden according to $p_{n-1,k}$. 
Hence,
\begin{align*}
W([n],k) & = \frac{c_n S_{k-1}([n-1])}{S_k([n])} W([n-1],k-1) + \frac{S_k([n-1])} {S_k([n])} ( c_n +  W([n-1],k))  \\
&= \frac{c_n S_{k-1}([n-1])}{S_k([n])} \frac{(k-1) S_{k}([n-1])}{S_{k-1}([n-1])} + \frac{S_k([n-1])} {S_k([n])} \left( c_n +  \frac{k S_{k+1}([n-1])}{S_k([n-1])} \right)  \\
& \qquad (\mbox{by induction})&\\
& = \frac{ c_n (k-1) S_{k}([n-1]) + c_n S_k([n-1]) + k S_{k+1}([n-1])}{S_k([n])} \\
& = \frac{ k (c_n S_k([n-1]) + S_{k+1}([n-1]))}{S_k([n])} \\
& =  \frac{k S_{k+1}([n])}{S_k([n])},
\end{align*}
which completes the proof.

We note that the lemma could also be proved directly from Lemma 2.3 of~\cite{lidbetter}.
\hfill $\Box$

\end{document}